\documentclass{article}

\usepackage{amsmath}
\usepackage{amsfonts}
\usepackage{amssymb}
\usepackage{graphicx}
\usepackage{hyperref}
\usepackage{url}
\usepackage{mathrsfs}

\begin{document}

\title{A Quasi-Unary  Representation of  Discrete Taxicab Geometry}
\author{Shahid Nawaz\thanks{email: \texttt{snawaz@albany.edu}}
\\
{\small Department of Physics, University at Albany-SUNY, }\\
{\small Albany, NY 12222, USA.}
}
\date{}
\maketitle
\begin{abstract}
In this paper we represent $n-$dimensional discrete Taxicab geometry by base--($4n+1$) numeral system. The algebraic structure of this base--($4n+1$) system is similar to unary system, we call it quasi-unary (QU) representation. QU representation generalizes translation and rotation to transform any geometrical object (shape) into another shape.
\end{abstract}
\section{Introduction}\label{Intro}
In the early 20th century  Herman Minkowski (1864-1909) published a whole family of metrics (the distance formula) defined on  Euclidean space. If $A=(x_1,x_2,x_3,\ldots, x_n)$ and $B=(x^\prime_1,x^\prime_2,x^\prime_3,\ldots, x^\prime_n)\in R^n$ are two points, then
\begin{equation}
d_k(A,B)=\left(\sum_{i=1}^{n}|x^\prime_i-x_i|^k\right)^{1/k}\,,
\end{equation}
 is the distance between  the two points. Here $k=1$ gives Taxicab metric and $k=2$ gives Euclidean metric. In this paper we are concerned with Taxicab geometry ($k=1$ case). Taxicab geometry is based on an ideal city in which all the streets are assumed to run either north and south or east and west \cite{Krause1986}. In a continuous taxicab geometry the blocks (buildings) are assumed to be of point size. In this paper we are particularly interested in discrete taxicab geometry (DTG) in which the blocks have  small but finite size.   Our main goal  is to represent the streets by strings of integers whose algebraic structure is similar to unary numeral system.  In standard numeral system one deals with the strings of symbol  1 (base-1 system). In this paper, on the other hand, we use more than one symbol. Therefore, we shall call it quasi-unary  (QU) representation.  An earlier version of this work can be found in \cite{SN2009}. The approach in this paper is, however, different.

The paper is organized as follows. The formalism is presented in section~\ref{formalism}. Length, inner-product, translation, rotation, shape transformation and functions are discussed in the subsequent sections. The paper is concluded in section~\ref{conclusions}.

\section{Formalism}\label{formalism}
\begin{figure}[h!]
  \centering
    \includegraphics[width=1.0\textwidth]{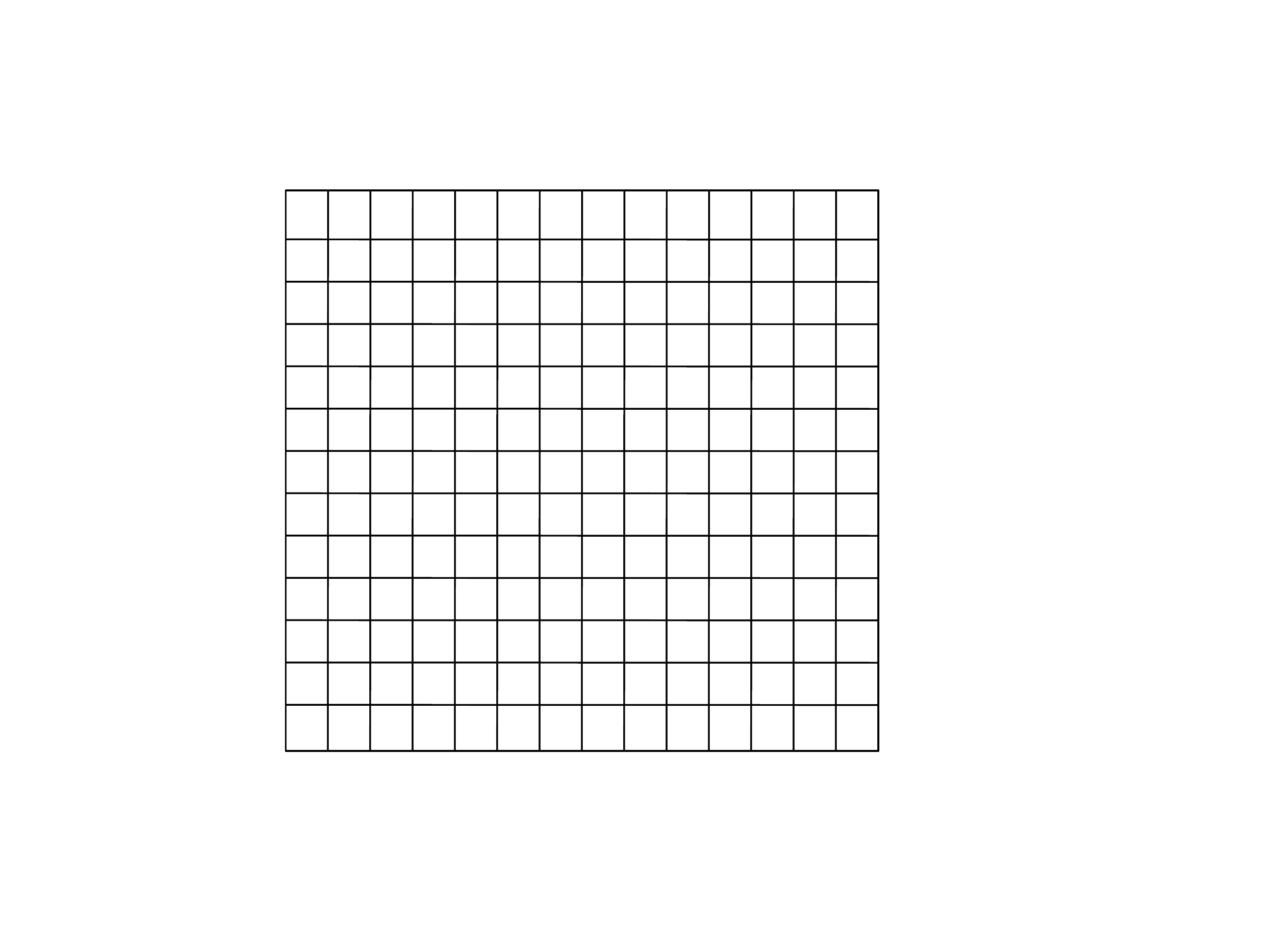}
       \vspace{-0.75in}
     \caption{A 2$-$dim discrete Euclidean space.}\label{grid}
\end{figure}

Consider an $n-$dimensional discrete Euclidean space (DES). A 2$-$dim  DES is shown in Fig.~\eqref{grid}. It consists of  pixels (blocks) and edges (streets of width zero). A taxicab can only run along the streets and we assume that it  can only pick up and drop passengers at the corners of a block. For example, if a cab goes 3 blocks to the east and then 4 blocks to the north then the total distance covered by the cab is 7 block. This is the taxicab distance. In contrast, the Euclidean distance would be 5 blocks.

Next we come to label the streets. We shall denote $x_1-$axis by $1-$axis, $x_2-$axis by $2-$axis, $x_3-$axis by $3-$axis, and $x_n-$axis by $n-$axis. The origin of  coordinate system is labeled as a dot `$\cdot$'.  The origin can be thought of as the reference point. A unit \emph{move} or a \emph{walk} parallel to   positive 1-axis is labeled as $1^+$. It describes a directed line from corner to corner of a block. An arrow is placed at the head of the line, that is, 
\begin{equation}
\rightarrow\ =1^+
\end{equation}
If a cab starts from the origin and moves one unit towards  the right,  it will be labeled as a dot followed by $ 1^+$,  
\begin{equation}
\cdot\!\!\rightarrow\ =. 1^+
\end{equation}
Similarly two consecutive moves towards the right can be obtained by concatenating $1^+$ with itself,
\begin{equation}
\cdot\!\!\rightarrow\!\rightarrow\ =.1^+1^+
\end{equation}
Our labeling scheme is analogous to \emph{directed graph}. However, there is a difference.  Only the edges are labeled. Except for the origin, the  vertices (corners) are left unlabeled. It is not, however. necessary that a walk  always starts from the origin. For example, $1^+.1^+$ is also possible which means that the initial point is one unit to the left of the origin (negative $1-$axis). In this case, start from the left go one unit towards the origin and then pass the origin to obtain $1^+\cdot1^+$. Thus the origin (dot) can float anywhere in the label set.

The moves/walks parallel to  negative directions and the other axes can be similarly defined. A unit walk parallel to  negative $1-$axis is labeled $1^-$.  Similarly $2^+$ ($2^-$) is a unit walk parallel to positive (negative) $2-$axis, and so on. It should be noted that $1^+$, $2^+$, \ldots work like \emph{basis}. The opposite of $1^+$ is $1^-$, for instance. They nullify the effect of each other when concatenated in a row,
\begin{equation}
1^+1^-=1^-1^+=0\,,
\end{equation}
where $0$ represents no move. This means that if $a$ and $b$ are strings, then
\begin{equation}
ab=a0b=a00b=\ldots
\end{equation}
where $ab$ is not a product but a concatenation of $a$ with $b$. Further, a string of finite $1^+$'s is labeled as $\overset{l}{\overline{ 1^+}}$ ( $l$ is the number of $1^+$ in the sequence). An infinite string of $1^+$'s is labeled as $\overset{\infty}{\overline{ 1^+}}\,\overset{\text{def.}}{=}\overline{ 1^+}$ .

\begin{figure}[h!]
\vspace{-10pt}
  \centering
    \includegraphics[width=1.0\textwidth]{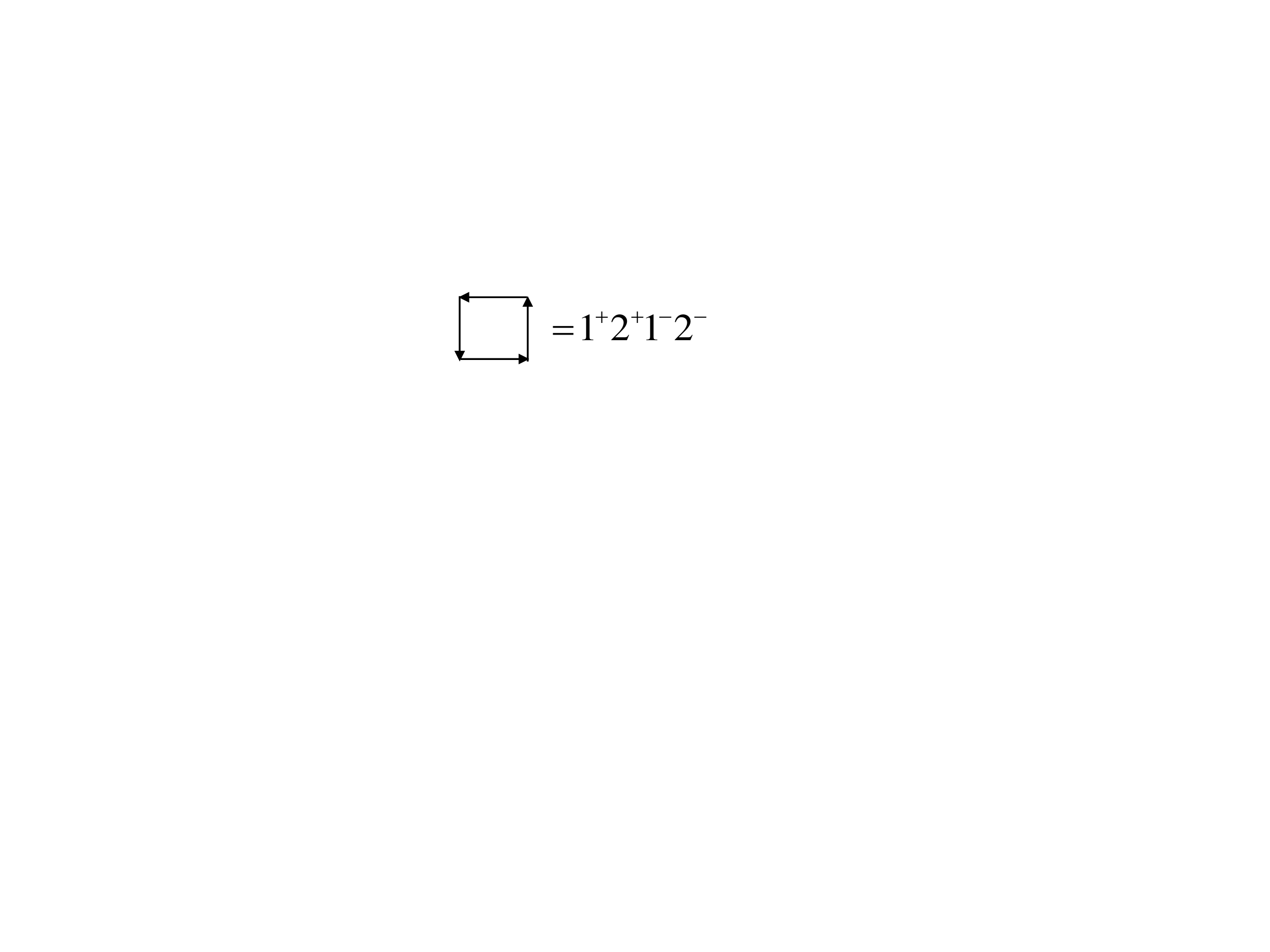}
       \vspace{-2in}
     \caption{A square}\label{square}
\end{figure}

Using the above method, we can draw any geometrical object (shape). A square, for example, is drawn in  Fig.~\eqref{square}. As we can see our labeling scheme involves concatenation of strings. It is similar to unary numeral system. But  standard  unary numbers involve strings of symbol 1 (base--1 numeral system). Our case, on the other hand, involves more than one symbol. Take, for example, the square in Fig.~\eqref{square} which involves four different symbols. We shall call our representation as quasi-unary (QU) representation. In fact QU representation is a base--($4n+1$) numeral system (see below).

Our labeling scheme so far only allows  us to draw connected geometrical objects. We can also draw disconnected geometrical objects that are made up of several pieces. We can do so by introducing new labels for missing links. Consider $i-$axis ($i=1,2,3,\ldots,n$). Let $\overset{k}{\overline{ i^+}}$ and $\overset{l}{\overline{ i^+}}$ be two line segments separated by a unit distance. The missing link or street or blank character, in the forward direction, is represented by $i^+_0$ ($i^+_0=1_0^+,2_0^+,\ldots,n^+_0$). Therefore the union of  $\overset{k}{\overline{ i^+}}$ and $\overset{l}{\overline{ i^+}}$ can be written as  $\overset{k}{\overline{ i^+}}\,i_0^+\,\overset{l}{\overline{ i^+}}$. This means that the two line segments are separated by $i^+_0$. Similarly $i^-_0$ can be understood. The symbols $i^\pm_0$ can be used for several purposes. They are also used for translation (see section \ref{translation}). They are  used also for leaving the axis by jumping off and then a unit distance away returning back to the axis.

To summarize this section, we note that any geometrical object can be labeled according to the map  $\alpha: \mathcal{G}\to \mathcal{M}$, where $\mathcal{G}$ is the set of geometrical objects such as line segments, curves etc., and $\mathcal{M}$ is the set of the corresponding labels. The elements of $\mathcal{M}$ will be  called strings. 
 
The set $\mathcal{M}$ has the following properties. Let $a,b,c\in\mathcal{M}$, then
\begin{enumerate}
\item \textbf{Additivity}
\begin{equation}
a\oplus b=ab\in\mathcal{M}\,,
\end{equation}
where $\oplus$ means concatenation. Here string $a$ is concatenated from left with string $b$ to obtain  string $ab$. String concatenation is non-commutative ($ab\neq ba$) in general. 
\item \textbf{Associativity} 
\begin{equation}
a\oplus(b\oplus c)=(a\oplus b)\oplus c\,.
\end{equation}
\emph{Proof:}
\begin{eqnarray}
l.h.s.&=&a\oplus(b\oplus c)\notag\\
&=&a\oplus bc\notag\\
&=&abc\notag
\end{eqnarray}
\begin{eqnarray}
r.h.s.&=&(a\oplus b)\oplus c\notag\\
&=&ab\oplus c\notag\\
&=&abc\,.
\end{eqnarray}
\item \textbf{Multiplication by a scalar number}
\begin{equation}
l\ast ab=\overset{l}{\overline{a}}\,\overset{l}{\overline{b}}\,,
\end{equation}
where $l$ is a positive integer.
\item \textbf{Subtraction (removal)} 
\begin{equation}
ab\ominus b=a\,,
\end{equation}
and
\begin{equation}
\ominus a\oplus ab=b\,.
\end{equation}
\item \textbf{Additive inverse}
\begin{equation}
a\ominus a=0
\end{equation}
\item  \textbf{Additive identity}
\begin{equation}
0ab=ab=a0b=a00b=\ldots=ab0
\end{equation}
\item \textbf{Atomic strings}\\
 The set 
\begin{equation}
\mathcal{A}=\{1^\pm,2^\pm,\ldots, n^\pm\}\subseteq\mathcal{M}\,,
\end{equation}
 constitutes of atomic strings. Atomic strings can also be called as orthonormal basis. They satisfy the following property,
\begin{equation}
i^+i^-=i^-i^+=0,\quad i^\pm\in \mathcal{A}\,.
\end{equation}
\item \textbf{Base--(4n+1) numeral system}\\
 Quasi-unary representation is in fact a base--($4n+1$) system. The \emph{digits} of this system are given by
\begin{equation}
\mathcal{B}=\{0,1^-,1^+,1^+_0,1^-_0,2^-,2^+,2^+_0,2^-_0,\ldots,n^-,n^+,n^+_0,n^-_0\}\subseteq\mathcal{M}\,.\label{base 4nplus1}
\end{equation}
It contains atomic strings, blank characters and zero. Various combination of the elements of  $\mathcal{B}$ correspond to various shapes. 
\end{enumerate}

\section{Length}\label{length}
Length can be obtained according to the map $|.|:\mathcal{M}\to R$. Let $i-$axis be any axis, then the length of a line segment is given by
\begin{equation}
|\overset{l}{\overline{i^+}}|=|\overset{l}{\overline{i^-}}|=ls\,,
\end{equation}
where $s$ is the length of the atomic string and $ls$ is the total length of the line segment. Here $ls$ is a product but not a concatenation. Concatenation is only used for strings. 

The map $|.|:\mathcal{M}\to R$ has the following properties:
\begin{enumerate}
\item 
\begin{equation}
|ab|=|ba|, \quad a,b\in\mathcal{M}.
\end{equation}
\item
\begin{eqnarray}
|\overset{k_i}{\overline{ i^+}}\,\overset{l_i}{\overline{ i^-}}\,\overset{k_j}{\overline{ j^+}}\,\overset{l_j}{\overline{ j^-}}|=|k_i-l_i|s+|k_j-l_j|s\,,\ \ \text{for}\ i\neq j
\end{eqnarray}
where $i$'s and $j$'s are atomic strings. Here various axes are  weighted equally. Otherwise $|i|=s_i$ and $|j|=s_j$ for $i\bot j$. 
\item \textbf{Taxicab length}\\
 The taxicab length is given by,
\begin{eqnarray}
d_T&=&|\overset{k_1}{\overline{1^-}}\,\overset{l_1}{\overline{1^+}}|+|\overset{k_2}{\overline{2^-}}\,\overset{l_2}{\overline{2^+}}|+\ldots+|\overset{k_n}{\overline{n^-}}\,\overset{l_n}{\overline{n^+}}|\notag\\
&=&\sum_{i=1}^n|l_i-k_i|s
\end{eqnarray}
\end{enumerate}

\section{Inner-product} Let $i^\pm$ and $j^\pm$ be atomic strings (orthonormal basis), then the inner-product is given by,
\begin{equation}
\langle i^+, j^+\rangle=s^2\delta_{ij}
\end{equation}
\begin{equation}
\langle i^+, j^-\rangle=-s^2\delta_{ij}\,,
\end{equation}
\begin{equation}
\langle i^-, j^+\rangle=-s^2\delta_{ij}\,,
\end{equation}
\begin{equation}
\langle i^-, j^-\rangle=s^2\delta_{ij}\,,
\end{equation}

\section{Translation}\label{translation}
Let $.a$ be a string  whose one end lies on the origin, then under translation 
\begin{equation}
.a\to a^\prime=.\,\overset{l}{\overline{i_0^+}}\,a\,,
\end{equation}
where $\overset{l}{\overline{i_0^+}}$ is the gap between the origin and string $a$. This means that $a$ is moved $l$ units in the forward direction. 

Length is invariant under translation
\begin{equation}
|a|=|\overset{l}{\overline{i_0^+}}\,a|\,.
\end{equation}
Note that 
\begin{equation}
|\overset{l}{\overline{i_0^+}}|=ls\,,
\end{equation}
is the length of the gap (blank characters). But $|i^+_0a|$ is equivalent to $|0a|=|a|$.

\section{Rotation}\label{rotation}
Length is not invariant under rotation in taxicab geometry \cite{OEB2002}. In Taxicab geometry, length is only preserved by a rotation of $\pi/2,\pi,3\pi/2$ and $2\pi$. In taxicab geometry, straight lines are those which are parallel to a coordinate axis. When a straight line is rotated, it deforms to a steps-type line. 

In quasi-unary (QU) representation, rotation is very simple. Here we consider rotation of a line segment which is  initially parallel to a coordinate axis. Let $i,j,k$ and $l$ be atomic strings. Consider,
\begin{equation}
a=.\,\overset{p}{\overline{k}}\,,
\end{equation}
where $\overset{p}{\overline{k}}$ is a line segment whose one end lies on the origin and is parallel to $k-$axis. Rotating $a$ about $l-$axis, then under rotation
\begin{equation}
a\to a^\prime=.\,\overset{p/(q+r)}{\overline{\left
(\overset{q}{\overline{i}}\,\overset{r}{\overline{j}}\right)}}\,,
\end{equation}
that is
\begin{equation}
.\,\overset{p}{\overline{k}}\to.\,\overset{p/(q+r)}{\overline{\left
(\overset{q}{\overline{i}}\,\overset{r}{\overline{j}}\right)}}\,,\label{rotation}
\end{equation}
where $p,q,r$ are non-negative integers such that $p\geq q+r$, and $p/(q+r)$ is also an integer.\\

\textbf{Examples:} 
\begin{enumerate}
\item Rotation by $45^{\text{o}}$:  Let $k=1^+$, $i=1^+$, $j=2^+$, and $q=r=1$, then
\begin{equation}
R_{45^{\text{o}}}\left(.\,\overset{2p}{\overline{1^+}}\right)=.\,\overset{p}{\overline{1^+2^+}}
\end{equation}
Here the line segment $.\,\overset{2p}{\overline{1^+}}$ is rotated counterclockwise to obtained $.\,\overset{p}{\overline{1^+2^+}}$ Note that $.\,\overset{p}{\overline{1^+2^+}}$ is a steps-type line segment at $45^{\text{o}}$. In the limit when the step size approaches zero then $p\to\infty$ and so $\overset{p}{\overline{1^+2^+}}$ tends to the straight line $\overline{1^+2^+}$ .
\item Rotation by $90^{\text{o}}$: For $k=1^+$, $i=2^+$, $q=1$, and $r=0$, then
\begin{equation}
R_{90^{\text{o}}}\left(.\,\overset{p}{\overline{1^+}}\right)=.\,\overset{p}{\overline{2^+}}
\end{equation}
\item Rotation by $135^{\text{o}}$: For $k=1^+$, $i=1^-$, $j=2^+$, and $q=r=1$, then
\begin{equation}
R_{135^{\text{o}}}\left(.\,\overset{2p}{\overline{1^+}}\right)=.\,\overset{p}{\overline{1^-2^+}}
\end{equation}
etc.
\end{enumerate}

\section{Shape Transformation}\label{shape transformation}
Quasi-unary (QU) representation generalizes translation and rotation to transform any geometrical shape into another shape. Various shapes can be obtained from various combination of the elements of the set of digits. Recall eq.~\eqref{base 4nplus1},
\begin{equation}
\mathcal{B}=\{0,1^-,1^+,1^+_0,1^-_0,2^-,2^+,2^+_0,2^-_0,\ldots,n^-,n^+,n^+_0,n^-_0\}\subseteq\mathcal{M}\,,\tag{\ref{base 4nplus1}}
\end{equation} 
which is the of set of digits of base--($4n+1$) numeral system. Here we use base--($4n+1$)  system to transform a line segment into any other shape. The shape transformation (ST) is then given by,
\begin{equation}
\overset{p}{\overline{c_k}}\,\stackrel{ST}{\rightarrow}\,\overset{p/\sum_mq_m}{\overline{\left
(\overset{q_1}{\overline{c_1}}\,\overset{q_2}{\overline{c_2}}\ldots \overset{q_l}{\overline{c_l}}\right)}}\,,\quad c_k,c_m\in \mathcal{B}\label{ST}
\end{equation}
This generalizes eq.~\eqref{rotation}. Here $p,q_m$ are non-negative integers, $p\geq\sum_mq_m$, and $p/\sum_mq_m$ is also an integer.\\

\textbf{Example:} Consider the transformation of a line segment into two squares. Here $c_k=1^+,\  p=9,\  c_1=2^+,\  c_2=1^-,\  c_3=2^-,\ c_4=1^+,\ c_5=1^+_0,\ c_6=1^+,\ c_7=2^+,\ c_8=1^-,\ c_9=2^-$, and $q_1=q_2=\ldots=q_9=1$, then
\begin{equation}
\overset{9}{\overline{1^+}}\,\stackrel{ST}{\rightarrow}\,2^+1^-2^-1^+1^+_01^+2^+1^-2^-\,,\label{two squares}
\end{equation}
where $2^+1^-2^-1^+$ and $1^+2^+1^-2^-$ are two squares separated by $1^+_0$ (figure not shown).

\section{Functions}\label{function}
Quasi-Unary (QU) representation can also be extended to functions. Treat $x$-axis as $1-$axis, and represent $f(x)$ by  $2-$axis, then the graph of $f(x)$ is given by,
\begin{equation}
\{(x,y)|\,y=f(x)\}=\overset{p_1}{\overline{i}}\overset{q_1}{\overline{j}}\overset{p_2}{\overline{i}}\overset{q_2}{\overline{j}}\ldots\,,
\end{equation}
where  $i=1^\pm$, and $j=2^\pm$.\\

\textbf{Examples:}
\begin{enumerate}
\item 
\begin{equation}
\left\{(x,y)|\,y=\frac{m}{n}x;\, m,n\in N\,\,\text{and}\,\,x\geq0\right\}=.\,\overline{\overset{n}{\overline{1^+}}\,\,\overset{m}{\overline{2^+}}}\label{S Line}
\end{equation}
This means that when $x=n$, then  $f(x)=m$. Note that eq.~\eqref{S Line} describes a straight line whose one end is on the origin and lies in the first quadrant.
\item 
\begin{equation}
\left\{(x,y)|\,y=x^2;\ x\geq0\right\}=.\,1^+2^+1^+1^+2^+2^+2^+2^+\ldots,
\end{equation}
etc.

\end{enumerate}

\section{Conclusions}\label{conclusions}
In this paper we have shown that $n-$dimensional Taxicab geometry can be represented by base--($4n+1$) numeral system, we call it quasi-unary (QU) representation. Any element (number) of this base--$(4n+1)$ system corresponds to an extended object such as line segment, curve etc. Length, translation, rotation and functions are expressed in terms QU representation. We have shown that new objects (shapes) can be obtained by transforming old shape using shape transformation. In QU representation, shape transformation is the generalization of translation and rotation.

%
 

              \bibliographystyle{utphys}

\providecommand{\href}[2]{#2}\begingroup\raggedright\endgroup

\end{document}